\documentclass[graybox]{svmult}

\usepackage{type1cm}        
\usepackage{makeidx}         
\usepackage{graphicx}        
\usepackage{multicol}        
\usepackage[bottom]{footmisc}

\usepackage{algorithm2e}
\usepackage{tikz}

\usepackage{newtxtext}       %
\usepackage[varvw]{newtxmath}       

\usepackage[version=4]{mhchem}

\usepackage{accents}

\makeindex             

\newcommand{\rev}[1]{#1} 

\newcommand{\uh}{\ensuremath{\hat{\mathbf{u}}_h}}
\newcommand{\Muh}{\ensuremath{\mathcal{M}\left(\hat{\mathbf{u}}_h}\right)}
\newcommand{\Uh}{\ensuremath{\hat{\mathbf{U}}_h}}

\newcommand{\HU}{\ensuremath{H\left(\mathbf{U}\right)}}

\newcommand{\Ue}{\mathbf{U}} 
 
\newcommand{\ue}{\mathbf{u}} 
\newcommand{\eps}{\varepsilon}
\newcommand{\Rc}{\mathcal{R}_c}
\newcommand{\Rs}{\mathcal{R}_s}
\newcommand{\Reps}{\mathcal{R}_{\epsilon}}
\newcommand{\Repsc}{\mathcal{R}_{{\epsilon},c}}
\newcommand{\Rd}{\mathcal{R}_{\delta}}
\newcommand{\rs}{r_s}
\newcommand{\F}{\mathbf{F}}

\newcommand{\g}{\mathbf{g}}
\newcommand{\R}{\mathbb{R}}

\newcommand{\source}{\mathbf{R}}
\newcommand{\RN}{\mathbb{R}^M}
\newcommand{\Rn}{\mathbb{R}^m}
\newcommand{\proj}{\mathbb{P}}
\newcommand*\diff{\mathop{}\!\mathrm{d}}

\newcommand{\Linfstc}{\ensuremath{L^\infty}}
\newcommand{\Linfsts}{\ensuremath{L^\infty}}

\DeclareMathOperator{\jacobian}{D}
\DeclareMathOperator{\jacobianu}{D_{u}}
\DeclareMathOperator{\hessian}{D^2}

\DeclareMathOperator{\Dx}{\mathrm{d}x}
\DeclareMathOperator{\Dt}{\mathrm{d}\tau}

\DeclareMathOperator{\Id}{I}
\newcommand{\G}{\ensuremath{\mathcal{G}}}
\renewcommand{\I}{\ensuremath{\mathcal{I}}}
\newcommand{\V}{\ensuremath{\mathcal{V}}}

\newcommand{\nt}{n}

\newcommand{\ubar}[1]{\underaccent{\bar}{#1}}

\newcommand\Tstrut{\rule{0pt}{2.6ex}}         
\newcommand\Bstrut{\rule[-1.5ex]{0pt}{0pt}}   

\begin{document}

\title*{Model adaptation for hyperbolic balance laws}
\author{Jan Giesselmann, Hrishikesh Joshi, Siegfried Müller and Aleksey Sikstel*}
\institute{Jan Giesselmann, Hrishikesh Joshi \at \rev{Technical University of Darmstadt, Department of Mathematics}, Dolivostr. 15, 64293 Darmstadt,
Germany, \email{ \{giesselmann / joshi \}@mathematik.tu-darmstadt.de}
\and
Siegfried Müller \at RWTH Aachen University, Institut für Geometrie und praktische Mathematik,   Templergraben 55, 52056 Aachen,
Germany, \email{mueller@igpm.rwth-aachen.de}
\and
Aleksey Sikstel* corresponding author \at University of Cologne, Department of Mathematics, Weyertal 86 – 90, 50931 Köln, Germany, \email{a.sikstel@uni-koeln.de}}
%
%
\maketitle

\abstract{  In this work, we devise a  model
  adaptation strategy for a class
  of model hierarchies consisting of two levels of model complexity. In particular, the fine
model consists of a system of hyperbolic balance
laws with stiff reaction terms and the coarse model
consists of a system of hyperbolic conservation laws.
We employ the relative entropy stability framework  to obtain an a posteriori  modeling error
estimator.
The efficiency of the model adaptation strategy
is demonstrated by conducting simulations for chemically reacting fluid mixtures in one
space dimension.}
\vspace{-0.5cm}
\section{Introduction}\vspace{-0.3cm}
\label{sec:intro}

Simulating hyperbolic balance laws featuring stiff, non-linear source terms can be computationally expensive due to the small time step sizes required for \rev{the stability of explicit time stepping methods}, cf. \cite{CockburnShu98},~\cite{ssp1}, or the iterative nature of  implicit time stepping methods, cf.~\cite{russonew},~\cite{russo2},~\cite{russo}. In some cases, the system of equations can be simplified given some constraints hold, leading to a system of conservation laws. This gives rise to a model hierarchy consisting of two levels of complexity; a system of hyperbolic balance laws and a system of hyperbolic conservation laws.

 We propose a model \rev{adaptation} strategy based on a posteriori error analysis that relies on the relative entropy stability framework, cf.~\cite{Dafermos},~\cite{Tzavaras2005}. This stability framework requires one of the solutions being compared to be Lipschitz continuous~\cite{Dafermos}. Since the numerical solution will not generally have the necessary regularity, a reconstruction of the numerical solution needs to be computed. The reconstruction method proposed in~\cite{GiesselmannDedner2016}, in the context of a posteriori error analysis of discontinuous Galerkin methods for hyperbolic conservation laws, is employed. \rev{A similar approach was used by Giesselmann et.~al.~in~\cite{GiesselmannPryer2017} where a model adaptation strategy} employing the relative entropy stability framework for \rev{a} model hierarchy consisting of the Euler equations and  Navier-Stokes Fourier equations \rev{has been proposed}.


\rev{For simplicity of exposition we restrict ourselves to one spatial dimension.} Let the \textbf{complex system} be defined on the \rev{one}-dimensional torus $\mathbb{T}$  for the unknowns $\Ue: \mathbb{T}\times\mathbb{R}^+\rightarrow\mathbb{R}^M$ as follows:
\begin{align}\label{eq01}
  \begin{split}
    &\partial_t \Ue+\partial_x \mathbf{F}(\Ue)=\frac{1}{\eps}\mathbf{R}(\Ue),\quad\Ue(x,\, 0) = \Ue_0(x)
  \end{split}
\end{align}
\noindent
where $\eps> 0$, $\mathbf{F} \in \rev{C^2(\mathbb{R}^M,\mathbb{R}^M)}$ is the flux function and $\mathbf{R} \in \rev{C^2(\mathbb{R}^M,\mathbb{R}^M)}$ the source term.
Moreover, the complex system is assumed to be equipped with a strictly convex entropy $H:\mathbb{R}^M \rightarrow \mathbb{R}$ and an entropy flux $Q: \mathbb{R}^M \rightarrow \mathbb{R}$ satisfying the compatibility condition $\jacobian Q(\Ue) =  \jacobian \HU\cdot \jacobian \mathbf{F}(\Ue), ~~ \forall~ \Ue \in \R^M$. \rev{Thus, smooth solutions of~\eqref{eq01} satisfy the  entropy balance law
$ \partial_t \HU + \partial_x Q (\Ue) = \frac{1}{\eps} \jacobian \HU \mathbf{R}(\Ue)$ and we assume that the complex system satisfies the second law of thermodynamics, i.e. 
 $\frac{1}{\eps} \jacobian \HU\cdot \mathbf{R}(\Ue) \leq 0$.}

 At equilibrium, the source term vanishes and the model can be simplified. The governing equations of the \textbf{simple system} are obtained as follows. Projecting the complex system by a constant matrix $\mathbb{P} \in \mathbb{R}^{m\times M}$ such that $\mathbb{P}\source(\Ue)=\vec{0}$ with $\text{rank} ~ (\mathbb{P}) = m$ and defining $\mathbf{u}:=\mathbb{P}\Ue: \mathbb{T}\times\mathbb{R}^+\rightarrow \mathbb{R}^m$, we obtain $\partial_t \vec{u}+\partial_x \mathbb{P}\vec{F}(\vec{U})=0$. \rev{We assume that there exists a Maxwellian $\mathcal{M} \,\colon\,\mathbb{R}^m \rightarrow \mathbb{R}^M$ that  parameterizes} the equilibrium states by the conserved quantities $\ue \in \Rn$  as
 \vspace{-0.2cm}
\begin{equation}\label{S3}
  \mathbf{R}(\mathcal{M}(\ue)) = 0,~ \mathbb{P} \mathcal{M}(\ue) = \ue.\vspace{-0.2cm}
\end{equation}
\rev{At least formally}, in the limit $\eps\rightarrow 0$, the conserved quantities $\mathbf{u}$ satisfy the conservation law 
\vspace{-0.2cm}
\begin{equation}\label{eq03}
\partial_t \mathbf{u}+\partial_x \g (\ue)=0, \quad\ue(x,\, 0) = \ue_0(x)
\end{equation}
where the flux is given by $\g (\ue) := \mathbb{P}\mathbf{F} (\mathcal{M}(\ue)): \Rn \rightarrow \Rn$. Moreover, the gradient of the source term on the equilibrium manifold is assumed to satisfy the non-degeneracy condition $\text{dim} ~\text{Ker}\left(\jacobian \source (\mathcal{M}(\ue))\right) =m ,~ \text{dim} ~\text{Im}\left(\jacobian \source (\mathcal{M}(\ue))\right) =M -m$. This condition is required for the entropy consistency  of the entropy-entropy flux pair for the simple system, i.e.~the induced entropy and entropy pair $\eta (\mathbf{u}) := H(\mathcal{M}(\ue))$ and $q(\mathbf{u}):= Q(\mathcal{M}(\ue))$.  In~\cite{hrishikeshDiss} the entropy consistency has been shown by employing geometric properties of the projection $\mathbb{P}$ and the Maxwellian $\mathcal{M}$ with regard to the complex entropy $H$.

\rev{The structure of the model hierarchy we look at was introduced by Tzavaras in~\cite{Tzavaras2005} and extended in~\cite{MiroshnikovTrivisa2014}, where the relative entropy stability framework was employed to show convergence of solutions of the complex system to solutions of the simple system for  $\eps\to0$ as long as the simple system admits a sufficiently regular solution.
In this work we utilize the relative entropy framework to derive computable a-posteriori error estimates for this model hierarchy.}
\vspace{-0.8cm}
\section{Error estimator}\vspace{-0.3cm}
\label{sec:error-estimator}
 The relative entropy and entropy flux compares  states $\mathbf{U}, \mathbf{V} \in \RN$ as follows
 \begin{equation}\label{R1}
 H\left(\mathbf{U}|\mathbf{V}\right) := H (\mathbf{U}) - H (\mathbf{V}) - \jacobian H(\mathbf{V}) \cdot \left(\mathbf{U} - \mathbf{V}\right), 
 \end{equation}
 \begin{equation}\label{R2}
  Q\left(\mathbf{U}|\mathbf{V}\right) := Q (\mathbf{U}) - Q (\mathbf{V}) - \jacobian H(\mathbf{V}) \cdot \left(\F\left(\mathbf{U}\right) - \F\left(\mathbf{V}\right)\right).
\end{equation}
\rev{Then, the relative entropy dissipation is} defined as
\begin{equation}\label{E9}
   \mathfrak{D}(\Ue|\mathbf{V}) :=-\left(\jacobian \HU - \jacobian H(\mathbf{V}) \right)\cdot \left(\mathbf{R}(\Ue)-\mathbf{R}(\mathbf{V})\right)
  \end{equation}
  for states $\mathbf{U},\mathbf{V} \in \RN$. Furthermore, we assume that for every \rev{compact $\mathcal{B} \subset \RN$  there exists a $\nu\left(\mathcal{B}\right)>0$ such that for all $\mathbf{U},\mathcal{M}(\proj \Ue) \in \mathcal{B}$:}
   \begin{equation}\label{E12}
    \mathfrak{D} (\Ue|\mathcal{M}(\proj \Ue))\geq \nu |\Ue - \mathcal{M}(\proj \Ue)|^2.
   \end{equation}

Dynamic heterogeneous model adaptation involves decomposing the spatial domain $\mathbb{T}$ after each time step and employing the simple system in a sub-domain $\Omega_s$ and the complex system in $\Omega_c$ such that $\mathbb{T} = \overline{\Omega}_s \cup \overline{\Omega}_c$ and $\Omega_c \cap \Omega_s = \emptyset$. In addition, we define the set of interfaces by $\Gamma := \overline{\Omega}_s \cap \overline{\Omega}_c$. For the sake of clarity of the error estimator,  we fix the subdomains in what follows.


The relative entropy framework  requires the quantities being compared to the exact solution to be Lipschitz continuous. As, in general, the numerical solution \rev{does not have} the necessary regularity, we use the reconstruction method   introduced in~\cite{GiesselmannDedner2016} and~\cite{GiesselmannMakridakisPryer2015}. Thus, \rev{given numerical solutions $\vec{U}_h$ and $\vec{u}_h$ of the complex and the simple system, respectively, we denote their space-time Lipschitz continuous reconstructions on $\Omega_c$ and $\Omega_s$ by $\Uh: \Omega_c\times[0,T]\rightarrow \RN, ~~ \uh: \Omega_s\times[0,T]\rightarrow \Rn$.}

These reconstructions  satisfy the following perturbed systems of equations\vspace{-0.3cm}
\begin{equation}\label{RC1}
\partial_t \hat{\Ue}_h+\partial_x \mathbf{F}(\hat{\Ue}_h)-\frac{1}{\eps}\source(\hat{\Ue}_h)=:\Rc
\end{equation}
\begin{equation}\label{RC2}
\partial_t \uh+\partial_x \mathbb{P} \mathbf{F}(\mathcal{M}(\uh))=:\rs
\end{equation}
\begin{equation}\label{RC2b}
\partial_t \Muh+\partial_x \mathbf{F}(\mathcal{M}(\uh))=:\Rs
\end{equation}
with explicitly computable residuals $\Rc  \in L_2({\Omega_c}\times[0,T],\RN)$, $ \Rs \in L_2({\Omega_s}\times[0,T],\RN)$ and $\rs \in L_2({\Omega_s}\times[0,T],\Rn)$, where $\proj \Rs = \rs$.


The residual of the simple system, $\Rs$, contains information about both the discretization and modeling errors. However, we are interested in an error estimator that is able to distinguish between those. To this end, a residual splitting has been proposed in~\cite{hrishikeshDiss}:
\begin{equation}\label{res_decomp}
 \Rs = \Rd + \Reps, \quad\text{ where } \Reps: =  \Rs -  \jacobian \Muh \cdot \rs, ~~   \Rd :=  \jacobian \Muh \cdot \rs.
\end{equation}
$\Reps$ is referred to as the modeling residual and $\Rd$ as the weighted discretization residual of the simple system.

\begin{theorem}\label{err_est_thm}
  Let $\Ue\,\colon\, \mathbb{T}\times[0,T]\rightarrow  \RN $ with $T>0$  be \rev{an entropy} solution to the complex system~\eqref{eq01}. Let $\Uh$ and $\uh$ be Lipschitz continuous reconstructions of numerical solutions  on $\Omega_c$ and $\Omega_s$, respectively, as defined above. Furthermore, let $\Ue$, $\Uh$ and $M(\uh)$ take values in a convex set $\mathcal{D} \subset\R^M$  only and let\vspace{-0.25cm}
  \[
    \rev{\Uh(x,t) = \mathcal{M}(\uh(x,t)), \quad \forall x\in\Gamma,\, t\in [0,T]}\vspace{-0.25cm}
  \]
  hold at the interfaces. Then for any $0 \leq t \leq T$, there exist positive constants $C_{\ubar{H}},\, C_{\bar{H}},\, C_{\bar{F}}$ and $ C_{\bar{M}}$, independent of $\Ue$, $\Uh$ and $\uh$,  such that the following error estimate holds
\begin{align}\label{thm_error_est}
\int_{\Omega_c} & \left|\mathbf{U}(\cdot,t)-\Uh(\cdot,t)\right|^2 \Dx   +\int_{\Omega_s} \left|\mathbf{U}(\cdot,t)-\mathcal{M}(\uh(\cdot,t))\right|^2 \Dx  \\ \nonumber
& \leq \frac{1}{C_{\ubar{H}}}\left(I+D_c+D_s+ M_s\right)\exp\left(\frac{1}{C_{\ubar{H}}}\max\left(G_c, G_s\right)t\right),
\end{align}
where
\begin{align*}
  &I := \displaystyle{\int_{\Omega_s} H(\mathbf{U}|\mathcal{M}(\widehat{\mathbf{u}}_h))\Big|_{t=0} \Dx  + \int_{\Omega_c} H(\mathbf{U}|\hat{\mathbf{U}}_h)\Big|_{t=0} \Dx },  \\ 
   &D_c := \displaystyle{\frac{1}{2}\int_0^t \int_{\Omega_c}|\hessian H(\Uh)\cdot \Rc|^2\Dx  \Dt} ,\quad D_s:= \displaystyle {\frac{1}{2} \int_0^t\int_{\Omega_s}|\hessian \eta(\uh)\cdot \Rd|^2\Dx  \diff\tau},  \\
  &M_s := \displaystyle{\frac{\eps}{\nu}\int_0^t\int_{\Omega_s} |\hessian \eta(\uh)\cdot \Reps|^2\Dx \Dt  },  \\
  &G_c :=\displaystyle{\frac{1}{2}+ C_{\bar{F}}C_{\bar{H}} \|\partial_x\Uh\|_{\Linfstc} + C_{\bar{H}}\left( \Big\|  \frac{1}{\eps}\source(\Uh)\Big\|_{\Linfstc}+\Big\|  \frac{1}{\eps}\jacobian\source(\Uh)\Big\|_{\Linfstc}\right)},   \\ 
  &G_s: =\displaystyle{\frac{1}{2}+C_{\bar{F}} C_{\bar{H}} \|   \partial_x\Muh\|_{\Linfsts}+ C_{\bar{H}}  C_{\bar{M}} |\proj|^2  \|\Reps\|_{\Linfsts}}.
\end{align*}

\end{theorem}

\begin{remark}
  The constants $C_{\ubar{H}},\, C_{\bar{H}},\, C_{\bar{F}}$ and $ C_{\bar{M}}$ appearing in the error estimate originate from the following bounds employed in the proof of Theorem~\ref{err_est_thm}.
  Let $h_i: \RN\rightarrow\R$ and $m_k : \Rn \rightarrow \Rn$ be defined as
  \[
    h_i(\Ue) :=  \frac{\partial H(\Ue)}{\partial U_i},~~ i=0,\ldots,M, \quad\quad  m_k(\ue):= \jacobianu (\mathcal{M}(\ue))_k, ~~ k=0,\ldots,M.
  \]
  Then, by regularity assumptions there exist constants   $0 < C_{\bar{F}}, \, C_{\ubar{H}},\, C_{\bar{H}},\, C_{\bar{M}} < \infty$ such that  for all $\mathbf{V} \in \RN$ \vspace{-0.3cm}
  \begin{subequations}
    \begin{align}
      &\left(\sum_{i=0}^M\left(\mathbf{V}^T \hessian \F_i(\Ue) \mathbf{V}\right)^2 \right)^{\frac{1}{2}}\leq C_{\bar{F}} | \mathbf{V}|^2, \quad C_{\ubar{H}} | \mathbf{V}|^2 \leq  \mathbf{V}^T \hessian H(\Ue) \mathbf{V} \leq C_{\bar{H}} | \mathbf{V}|^2.\label{cons_H}\\
      &\left(\sum_{i=0}^M\left(\mathbf{V}^T\hessian h_i(\Ue)\mathbf{V}\right)^2\right)^{\frac{1}{2}} \leq C_{\bar{H}}| \mathbf{V}|^2,\quad\left(\sum_{k=0}^m\left(\left(\proj\mathbf{V}\right)^T\jacobianu m_k(\ue)\left(\proj\mathbf{V}\right)\right)^2 \right)^{\frac{1}{2}}\leq C_{\bar{M}}| \mathbf{V}|^2.\label{cons_m}
    \end{align}
  \end{subequations}
  The constants can be explicitly calculated given $\mathcal{D}, H, M$ and $\F$ and, thus, the right-hand side of the error estimate~\eqref{thm_error_est} is explicitly computable.
\end{remark}
In order to bound the distance between $\Ue$ and $\mathbf{U}_h$ we employ the  error splitting
\begin{align*}
|\mathbf{U}-\mathbf{U}_h| & \leq |\mathbf{U}-\Uh| +|\Uh - \mathbf{U}_h|, \\
|\mathbf{U}-\mathcal{M}\left(\mathbf{u}_h\right)| & \leq |\Ue-\mathcal{M}(\uh)| + |\mathcal{M}(\uh)-\mathcal{M}(\mathbf{u}_h)|.
\end{align*}
 The first terms on the right-hand side in the equations above can be bounded by the error estimates derived from the relative entropy framework, cf.~\cite{Dafermos},~\cite{Tzavaras2005}, and the second terms on the right-hand side are explicitly computable. For the complete proof of the theorem we refer to~\cite[Sec. 2.5.6]{hrishikeshDiss} in order to keep this work concise.

 \vspace{-0.7cm}
 \section{Model \rev{adaptation} strategy}\vspace{-0.3cm}
 \label{sec:model-adapt-strat}

Let the solution be approximated at times $0 = t^0 < t^1 <\ldots < T_F$  with an equidistant mesh with $\mathcal{N}_E$ number of cells on $\Omega = (a,b)$. The spatial cells  are given by
  $\V_i := (x_i,x_{i+1}),~~ i \in  \{1,\ldots,\mathcal{N}_E\}$,
where
$ x_i: = a + ih ,~ h := \frac{b-a}{\mathcal{N}_E}$. Moreover, let the numerical solution and its space-time Lipschitz-reconstruction in cells $\V_i\subset \Omega_c$ and $\V_j\subset\Omega_s$ be \rev{denoted by} $\mathbf{U}^{\nt}_i,\hat{\mathbf{U}}^{st,\nt}_i$ and $\mathbf{u}^{\nt}_j,\hat{\mathbf{u}}^{st,\nt}_j$, respectively. Then, the modeling error indicator for the simple system reads
\vspace{-0.3cm}
  \begin{equation}\label{ms_defn}
 \mathcal{\mathtt{M}}_s^{\nt,i}:= \left(\frac{1}{\Delta x_i}\cdot\frac{1}{t^{\nt+1}-t^\nt}\cdot\frac{\eps}{\nu}\int_{t^\nt}^{t^{\nt+1}}\int_{\V_i} |\hessian \eta(\hat{\mathbf{u}}^{st,\nt}_i)\cdot \Reps^i|^2\Dx \Dt\right)^{\frac12},  \vspace{-0.3cm}
\end{equation}
which is based on the $M_s$ term in the error estimate provided in Theorem~\ref{err_est_thm}, scaled by the mesh width and time step to account for the size of the space-time cell $\V_i\times [t^{\nt-1},t^\nt]$. This indicator allows to quantify the extent to which the simplifying assumption made to derive the simple system holds for cells in $\Omega_s$.

Let $\tilde{\mathbf{u}}_c^{\nt,i}$ be a numerical solution \rev{projected} from the complex to the simple system. Let $\hat{\mathbf{u}}^{st,\nt}_{c,i}:\V_i \times [t^{\nt-1},t^\nt] \rightarrow \Rn$ be a space-time Lipschitz continuous reconstruction of $\tilde{\mathbf{u}}_c^{\nt,i}$. Then the model coarsening residual is defined as\vspace{-0.1cm}
\begin{equation}\label{res_defc}\vspace{-0.1cm}
  \Repsc^i:=\left(\Id - \jacobian  \mathcal{M}\left(\hat{\mathbf{u}}^{st,\nt}_{c,i}\right) \cdot  \mathbb{P} \right)\cdot \partial_x \mathbf{F}\left(\mathcal{M}\left(\hat{\mathbf{u}}^{st,\nt}_{c,i}\right)\right)
\end{equation}
and the model coarsening error indicator reads\vspace{-0.3cm}
  \begin{equation}\label{mc_defn}\vspace{-0.3cm}
    \mathcal{\mathtt{M}}_c^{\nt,i}:= \left(\frac{1}{\Delta x_i}\cdot\frac{1}{t^{\nt+1}-t^\nt}\cdot\frac{\eps}{\nu}\int_{t^\nt}^{t^{\nt+1}}\int_{\V_i} |\hessian
      \eta(\hat{\mathbf{u}}^{st,\nt}_{c,i})\cdot \Repsc^i|^2\Dx \Dt \right)^{\frac12}.
\end{equation}
  $\mathcal{\mathtt{M}}_c^{\nt,i}$ is the complex system counterpart of $\mathcal{\mathtt{M}}_s^{\nt,i}$, i.e. $\mathcal{\mathtt{M}}_c^{\nt,i}$ allows to determine the extent to which the simplifying modeling assumption would hold in the case that we do decide to coarsen the model. If the simplifying assumption does not hold to the extent we prescribe, we continue employing the complex system.

\rev{When the model is coarsened, the numerical solution is converted from the complex to the simple system. Since we have employed the relative entropy framework for the error estimates, it suggests itself to also employ it for measuring coarsening distances.  Thus,  we define the model coarsening distance as}  \vspace{-0.3cm}
 \begin{equation}\label{keps_defn}\vspace{-0.3cm}
 \kappa_{\eps}^{\nt,i} := \left(\frac{1}{\Delta x_i}\cdot\frac{1}{t^{\nt}-t^{\nt-1}}\int_{\V_i} H\left(\mathbf{U}_h^{\nt,i}|\mathcal{M}(\proj \mathbf{U}_h^{\nt,i})\right) \Dx\right)^{\frac12}.
\end{equation}
We refer to~\cite[Sec. 3.4]{hrishikeshDiss} for a detailed discussion of the error indicators properties.
  
Model adaptation is done employing the model error indicators and the model coarsening distance. Since they provide two different pieces of information, one tolerance $\tau_r$ is set for the error indicators and another tolerance $\tau_\kappa$ for the coarsening distances.  Furthermore, when coarsening the model, we employ $\mathcal{\mathtt{M}}_c^{\nt,i}$ to ascertain if the simplifying assumption holds to the extent we prescribe. Switching from the complex to the simple system and back to the complex system in a matter of a few time steps should be avoided. This can occur when the modeling error indicators are close to the tolerance. To avoid switching frequently between the two systems, a factor of safety  $0< \mathfrak{f}_\eps< 1$ is employed when comparing $\mathcal{\mathtt{M}}_c^{\nt,i}$ to the tolerance.

\noindent
The models to be employed then are determined using Algorithm \ref{alg:model_adapt}.\vspace{-0.4cm}
\begin{algorithm}[!htb]
  \caption{Spatial model adaptation}\label{alg:model_adapt}
  \DontPrintSemicolon
  \KwIn{$\Theta_i^{\nt},\Theta_i^{\nt+1},\mathcal{\mathtt{M}}_{c/s}^{\nt,i},\kappa_\eps^{\nt,i}$}
   \lIf{$\Theta_i^{\nt}  = 1 ~~\text{and}~~ \mathcal{\mathtt{M}}_c^{\nt,i} < \mathfrak{f}_\eps\cdot\tau_r ~~\text{and}~~ \kappa_\eps^{\nt,i} < \tau_\kappa$}
   {
    \quad\quad$\Theta_i^{\nt+1} = 0$
   }
   \lElseIf{$\Theta_i^{\nt}  = 0 ~~\text{and}~~ \mathcal{\mathtt{M}}_s^{\nt,i} > \tau_r$}
   {
    \,\,\quad\quad\quad\quad\quad\quad\quad\quad$\Theta_i^{\nt+1} = 1$
   }
   \lElse
   {
    ~~~\,\quad\quad\quad\quad\quad\quad\quad\quad\quad\quad\quad\quad\quad\quad\quad\quad\quad\quad\quad\quad$\Theta_i^{\nt+1} = \Theta_i^{\nt}$
   }
   \Return{$\Theta_i^{\nt+1}$}
 \end{algorithm}\vspace{-0.5cm}
 Here we process a model indicator array $\Theta^\nt$ that stores entries
 \begin{equation*}
 \Theta^\nt_i :=  \begin{cases}
  1 , ~\text{if \rev{the} complex system is employed in}~ \V_i \in \G^\nt \\ 
  0 , ~\text{if \rev{the} simple system is employed in}~ \V_i \in \G^\nt
 \end{cases}.
\end{equation*}
 In the next step, the models to be employed are adjusted as follows:  The patch of cells where the simple system is \rev{used} cannot be a single cell since  there is a high likelihood for undesired frequent switching between the two systems that are triggered by the complex system set in the surrounding cells. 

The conversion of a numerical solution from the complex to the simple system is realized by projecting with the matrix $\mathbb{P}$. For the other direction it is required to evaluate the Maxwellian $\rev{\mathcal{M}}$, i.e.~solving a system of nonlinear equations, cf.~\cite[Sec. 3.3.4]{hrishikeshDiss}.
\vspace{-0.7cm}
\section{Numerical Results}\vspace{-0.3cm}
\label{sec:numer-exper}

We conduct numerical experiments for the  dissociation of oxygen, where the molecular oxygen dissociates into atomic oxygen with nitrogen acting as a catalyst:
\begin{center}\label{chem_reac}
\ce{$O_2 + N_2$ <=>$2O + N_2$ } .
\end{center}

\noindent
\textbf{Complex system.} Let the vector of  unknown quantities of the complex system be  $\Ue := \begin{bmatrix}\rho_{O_2},&\rho_O,&\rho_{N_2},&\rho v,& \rho E\end{bmatrix}^T$, where $\rho$ denotes the partial density of the \rev{species indicated by the subscript},   $\rho v$ the total momentum and  $\rho E$ the total energy. \rev{Furthermore, we enumerate the partial densities by $\rho_1 := \rho_{O_2}$, $\rho_2 : = \rho_O$ and $\rho_3 = \rho_{N_2}$}. The equations describing the local balance are
\begin{subequations}
  \label{eq:complexs-sys}
  \begin{align}
    &\partial_t \rho_k+\partial_x(\rho_k v)=m_k \sum_{j=1}^{\mathcal{N}_r}\nu^j_k \mathfrak{R}_j,\label{gov_rhoi}\\
    &\partial_t(\rho v)+\partial_x \left(\rho v^2 +p\right)=0, \label{gov_mom}\\
    &\partial_t(\rho E)+\partial_x \left( \left(p+\rho E\right)v \right)=0,\label{gov_tot_ener}
  \end{align}
\end{subequations}
\noindent
where $\nu_k^j=\beta_k^j-\alpha^j_k$ are constants with $\alpha_k^j, \beta_k^j \geq 0$  the stoichiometric coefficients of the reactions. $m_k$ is the molecular mass of $k$-th constituent and $p$ is the total pressure of the fluid mixture.  This system is a Class-I model according to Bothe and Dreyer~\cite{Bothe2015}.

Let $x_k := \frac{c_k}{c}$, where $c_k:=\frac{\rho_k}{m_k}$ and $c=\sum_{k=1}^{\rev{3}}c_k$, be the molar concentration of the $k$-th species. Then the reaction source term can be expressed as $ \mathfrak{R}_j = k^f_j \left(\prod_{i=1}^{\rev{3}} x_k^{\alpha_k^j} - \frac{1}{k^{eq}_j} \prod_{i=1}^{\rev{3}} x_k^{\beta_k^j}\right)$ with $k^f_j$ the forward and $k^{eq}_j$ the equilibrium reaction rate constants.

We assume the ideal gas law, i.e.~the partial pressures read $p_k = \frac{R}{m_k}\rho_kT$, where $R$ is the universal gas constant and $T$ the temperature. Furthermore, the partial entropy density is given by $s_k=s_k^R + c_{v,k}\ln\left( \frac{T}{T^R}\right)  - \frac{R}{m_k}\ln\left( \frac{\rho_k}{\rho_k^R} \right)$ where the superscript $^R$ denotes reference quantities and $c_{v,k}$ is the specific heat at constant volume. Thus, the total pressure reads $p = RT\sum_{k=1}^{\rev{3}} \frac{\rho_k}{m_k}$ and the entropy $\rho s = \sum_{k=1}^{\rev{3}}\rho_ks_k$. In~\cite{hrishikeshDiss} it was shown that $H=-\rho s$ and $Q=\rho sv$ is a valid (mathematical) entropy--entropy flux pair for the complex system~\eqref{eq:complexs-sys}.

\noindent
\textbf{Simple system.} The fluid mixture is in chemical equilibrium when the reaction terms of the complex system vanish.
The vector of unknown quantities of the simple system is given by  $\ue:=\begin{bmatrix}\rho_{O_2}+\rho_O,&\rho_{N_2},&\rho v,& \rho E\end{bmatrix}^T $ with the projection matrix given by
\begin{equation}\label{proj_d}
 \proj :=\begin{bmatrix}
          1 & 1 & 0 & 0 & 0 \\
          0 & 0 & 1 & 0 & 0 \\
          0 & 0 & 0 & 1 & 0 \\
          0 & 0 & 0 & 0 & 1
         \end{bmatrix}.
\end{equation}

\noindent
The forward and the equilibrium reaction rate constant is assumed to be of the form 
\begin{equation}\label{kf}
k_f := C \cdot T^{-2}\cdot \exp\left(\frac{-E}{T}\right), \quad k^{\text{eq}}_j=\exp\left(-\frac{1}{RT}\sum_{k=1}^{\mathcal{N}_c}\nu_k^j m_k g_k \right)
\end{equation}
with the constants $C = 2.9\cdot10^{13}~m^3mol^{-1}s^{-1}$  and $E = 597.5~K$. Furthermore, $g_k$ denotes the  Gibbs free energy $g_k(T,p):= e_{0,k} - Ts_k^R + c_{v,k}(T-T^R) - c_{p,k}\ln\left( \frac{m_k}{RT^R\rho_k^R}p  \right)$ where $e_{0,k}$ denotes the reference specific internal energy, $c_{p,k} := c_{v,k} + R$ the specific heat at constant pressure. Forward reaction rate constant of the form~\eqref{kf} is commonly employed in non-equilibrium hypersonic airflows, cf.~\cite{Nasa}. 

The mass of atomic oxygen and atomic nitrogen is \rev{$0.016 ~kg~mol^{-1}$} and \rev{$0.014~kg~mol^{-1}$}, respectively. The value of the specific gas constant is $R = 8.314 ~ J ~ K^{-1} ~ mol^{-1}$ and the reference temperature and pressure are assumed to be $T^R =2000K, p^R= 1.01325\times10^5 m^{-1}~kg\cdot s^{-2}$, respectively.

 The remaining thermodynamic and physical constants of the chemical constituents at hand are listed in Table \ref{table_cons}, cf.~\cite{thermochemical1} and~\cite{thermochemical2}.\vspace{-0.5cm}
 \begin{table}[h!]
\centering
\caption{Physical and thermodynamics constants}
\begin{tabular}[t]{|c|c|c|c|}
  \hline                                         $k$                            & $\rho_{O_2}$          & $\rho_O$              & $\rho_{N_2}$        \Tstrut\Bstrut  \\ \hline
  $\alpha_k$ [$kg$]                                                                    & 1              & 0             & 1            \Tstrut\Bstrut  \\ \hline
  $\beta_k$ [$kg$]                                                                    & 0              & 2              & 1             \Tstrut\Bstrut  \\ \hline
$m_k$ [$kg$]                                                                    & $0.032 $              & $0.016 $              & $0.028$             \Tstrut\Bstrut  \\ \hline
  $c_{v,k}$ [$J ~ mol^{-1}$]                                                        & $\frac{5}{2M_{O_2}}R$ & $\frac{3}{2M_{O_2}}R$ & $\frac{5}{2M_{O_2}}R$ \Tstrut\Bstrut\\ \hline
  $e_{0,k}$ [$J ~ K^{-1}~ mol^{-1}$]                                                & $249200$              & $0$                   & $0$                  \Tstrut\Bstrut \\ \hline
$\rho^R_k$ [kg]                                                                 & 1145                  & 1141                  & 1308                  \Tstrut\Bstrut\\ \hline
$s^R_k$ [$ J ~ K^{-1}~ mol^{-1}$]                                               & 205.15                                  & 161.1                 & 191.61                      \Tstrut\Bstrut\\ \hline
\end{tabular}
\label{table_cons}
\end{table}\vspace{-0.3cm}

\noindent
\textbf{Numerical scheme.} We employ a RKDG scheme, \rev{using a third order explicit SSP-RK time-stepping method and  quadratic polynomials for the spatial discretization. Furthermore, we use the local Lax-Friedrichs numerical flux and the minmod limiter introduced in~\cite{CockburnShu98}}. The solver is implemented in the \texttt{Multiwave} library cf.~\cite{MullerMultiwave2018}, ~\cite{GerhardIaconoMayMullerSchafer2015} and~\cite{gerhard2020waveletfree}.  \rev{The computational domain is $\Omega=\mathbb{T}$, i.e.~we set~periodic boundary conditions, and the number} of cells in the uniform mesh is  $\mathcal{N}_E = 1280$. The size of the time step is fixed throughout the computations and is set such that it satisfies the $CFL\leq 0.1$.  In the first time step, the complex system is \rev{used} in the entire domain $\Omega$. The model adaptation and the according domain decomposition is performed at the end of each time step using Algorithm~\ref{alg:model_adapt}. At the interfaces $\Gamma$ between the sub-domains $\Omega_c$ and $\Omega_s$ we employ the coupling strategy from~\cite[Sec. 3.3.6]{hrishikeshDiss}.

\noindent
\textbf{Shock tube test case.} The temperature and velocity is set to $T = 2000~K , v = 0~m/s$ in the entire computational domain. The pressure of the fluid mixture and the density of atomic oxygen is set to \vspace{-0.3cm}
\begin{equation}
  \begin{cases}
    p = 2\cdot 10^6,\quad \rho_0 = 0.01 & ~ \text{for}~ |x| \leq 0.5 \\
    p =  10^6,\quad \,\,\,\,\,\,\, \rho_0 = 0.005 & ~ \text{for}~ |x| \geq 0.5
  \end{cases}.
\end{equation}
Based on the values of the temperature, the pressure and the velocity, the equilibrium value $\mathbf{U}_{eq}$ is calculated.  The error indicator tolerance is set to $\tau_{r} = 0.16$, the coarsening distance tolerance to $ \tau_{\kappa} = 0.0016$ and the factor of safety is  $\mathfrak{f}_\eps = 0.25$.

\begin{figure}[!htb]
  \centering
  \includegraphics[width = 1\textwidth]{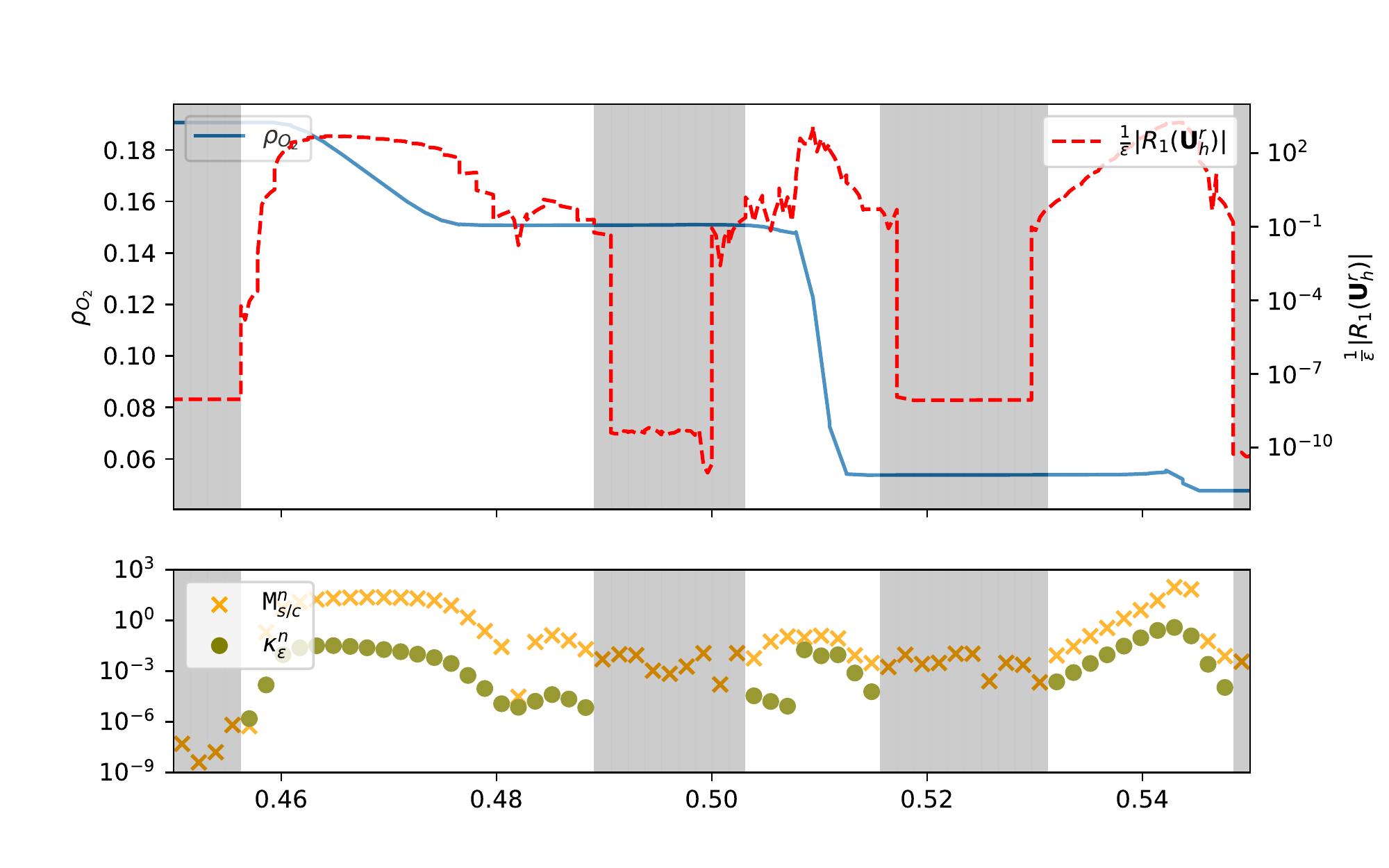}
  \caption{$t = 4.375 \cdot 10^{-4}$. Gray areas indicate $\Omega_s$, white background $\Omega_c$.}\label{fig:shock-tube}
\end{figure}

Figure~\ref{fig:shock-tube} depicts the model-adaptive numerical solution zoomed on the subdomain $[0.45, 0.55] \subset \Omega$ since outside of it the solution is constant.  In addition, the error indicators and the source term $\frac{1}{\varepsilon}|\vec{R}_1(\vec{U}^r_h)|$ of a reference simulation $\vec{U}^r_h$, where the complex system has been imposed in the whole domain $\Omega$, are presented.

The initial discontinuity gives rise to a rarefaction wave, a shock and a contact discontinuity.
 After the first time step, the model adaptation strategy leads to the simple system being employed in regions where the solution is constant. In regions where the waves emerge due to the initial discontinuity the complex system is selected, see Figure~\ref{fig:shock-tube}.

 As the contact discontinuity and the shock travels to the right and the rarefaction travels to the left, the numerical solution approaches equilibrium in the two plateaus between the three waves. As a result, the simple system is selected between the contact discontinuity and the shock from time $2.6125 \cdot 10^{-4}$ onward, and between the rarefaction and the contact discontinuity from time $t = 3.725 \cdot 10^{-4}$ onward. The switch to the simple system is made once the modeling error indicator and the model coarsening distance are smaller than the prescribed tolerances and the number of cells in the patch where the simple system is to be employed in more than one cell.

 The complex system has to be employed near the contact discontinuity, the shock and the rarefaction. In the case of the rarefaction, this is due to the fact that the rarefaction continues to expand as it travels to the left. Hence, the states that make up the rarefaction are away from the equilibrium manifold. In the case of the contact discontinuity and the shock, since the discontinuities  are spread over a few cells, the states of the discontinuities are away from the equilibrium manifold as well. The size of the two patches of cells in between the emerged waves where the simple system is employed increases in the long run. 

 

 Thus, the model adaptation strategy leads to employment of the simple system, locally in time and space and does not result in undesirable back and forth switching between the two systems. The switch back to the complex system is only triggered if the number of cells in a patch where the simple system is employed falls below one. We can infer that the model adaptation strategy works well by ensuring that the numerical solution is close to the equilibrium manifold and ensuring that the dynamics is close to that of the equilibrium dynamics before switching to the simple system. This is also \rev{suggested} by the fact that the source term in the reference simulations is of a smaller magnitude in the regions where the simple system is set.

 \bibliography{ref} 
\bibliographystyle{styles/spmpsci}
\end{document}